\documentclass[12pt]{article}

\usepackage{a4wide,epsfig}
\usepackage{amsfonts,amsmath}

\usepackage{pgfplots}
\pgfplotsset{compat=newest}

\newtheorem{theorem}{Theorem}[section]
\newtheorem{proposition}[theorem]{Proposition}

\newtheorem{lemma}[theorem]{Lemma}

\newcommand{\proof} [1]
   { \noindent {\bf Proof.} #1 \hfill\rule{0.5em}{1.2ex} \par\medskip}

\newcommand{\R}{\mathbb{R}}

\numberwithin{equation}{section} 

\begin{document}

\setcounter{page}{1}

\title{Towards coercive boundary element methods \\[1mm] for the wave equation}
\author{Olaf~Steinbach$^1$, Carolina~Urz\'ua--Torres$^2$, Marco~Zank$^3$}
\date{$^1$Institut f\"ur Angewandte Mathematik, TU Graz, \\ 
Steyrergasse 30, 8010 Graz, Austria \\[1mm]
{\tt o.steinbach@tugraz.at}  \\[3mm]
$^2$Delft Institute of Applied Mathematics, TU Delft, \\
Mekelweg 4, 2628CD Delft, The Netherlands \\[1mm]
{\tt C.A.UrzuaTorres@tudelft.nl} \\[3mm]
$^3$Fakult\"at f\"ur Mathematik, Universit\"at Wien, \\
Oskar-Morgenstern-Platz 1, 1090 Wien, Austria \\[1mm]
{\tt marco.zank@univie.ac.at}}

\maketitle

\begin{abstract}
  In this note, we discuss the ellipticity of the single layer boundary
  integral operator for the wave equation in one space dimension. This
  result not only generalizes the well-known ellipticity of the energetic
  boundary integral formulation in $L^2$, but it also turns out to be
  a particular case of a recent result on the inf-sup stability of boundary
  integral operators for the wave equation. Instead of the time derivative
  in the energetic formulation, we use a modified Hilbert transformation,
  which allows us to stay in Sobolev spaces of the same order. This
  results in the applicability of standard boundary element
  error estimates, which are confirmed by numerical results.
\end{abstract}

\section{Introduction}
Time-domain boundary integral equations and boundary element methods for
the wave equation are well established in the literature; we mention
the groundbraking works of Bamberger and Ha Duong \cite{Bamberger:1986},
Aimi et al. \cite{Aimi:2009}, and the review article
\cite{Costabel:2017} by Costabel and Sayas. Other works include
\cite{Gimperlein:2017, Hassel:2017, Joly:2017, Poelz:2021, Sayas:2013, Sayas:2016},
to mention a few.

The main difficulties in the numerical analysis of these formulations
are in the so-called norm gap, coming from continuity and coercivity
estimates in different space-time Sobolev norms. When using the
energetic boundary element method, a complete stability and error
analysis can be done in $L^2(\Sigma)$, see \cite{Joly:2017}, where
$\Sigma$ is the lateral boundary of the space-time domain
$Q := \Omega \times (0,T)$.

Using a generalized inf-sup stable variational formulation
\cite{SteinbachZank:2021a} for the wave equation,
in \cite{SteinbachUrzua:2021} we derived
inf-sup stability conditions for all boundary integral operators
in related trace spaces. In fact, this work was motivated by our
previous result \cite{Urzua:2021} on the spatially one-dimensional
case. When replacing the time derivative in the energetic
boundary integral formulation by a modified Hilbert transformation
\cite{SteinbachZank:2018}, the resulting composition with the
single layer boundary integral operator becomes elliptic in the
natural energy space $[H^{1/2}_{,0}(\Sigma)]'$, similarly to what is known
for boundary integral operators for second-order elliptic partial
differential equations. Note that $H^{1/2}_{,0}(\Sigma) :=
[H^1_{,0}(\Sigma),L^2(\Sigma)]_{1/2}$ is defined by interpolation,
with $H^1_{,0}(\Sigma) = \{ v \in H^1(\Sigma) : v(T)=0 \}$.
Analogously, $H^1_{0,}(\Sigma)$ covers zero initial conditions, i.e.,
$v(0)=0$.

In this paper, we present a detailed derivation of this new approach, and we
discuss the corresponding numerical analysis of a related new boundary
element method. In Section~2, we recall the energetic space-time boundary
integral formulation \cite{Aimi:2009,Joly:2017}, and we provide a
simplified proof of the ellipticity result in $L^2(\Sigma)$. 
In particular, we obtain that the single layer boundary integral
operator $V : L^2(\Sigma) \to H^1_{0,}(\Sigma)$ is an isomorphism.
Using duality arguments, we obtain that
$ V : [H^1_{,0}(\Sigma)]' \to L^2(\Sigma)$ is also an isomorphism.
Finally, by an interpolation argument, we conclude
that $ V:[H^{1/2}_{,0}(\Sigma)]' \to H^{1/2}_{0,}(\Sigma)$ is an
isomorphism as well. While this implies an inf-sup stability estimate,
as also discussed in \cite{SteinbachUrzua:2021}, in Section~3 we introduce a
modified Hilbert transformation ${\mathcal{H}}_T : H^{1/2}_{0,}(\Sigma) \to
H^{1/2}_{,0}(\Sigma)$, see \cite{SteinbachZank:2018}, 
to establish ellipticity of ${\mathcal{H}}_T V$ in
$[H^{1/2}_{,0}(\Sigma)]'$ in Section~4. Although the main result, as given in
Lemma \ref{Lemma HV}, still involves some unknown constant,
Proposition \ref{Proposition C} gives numerical evidence on the
behavior of the ellipticity constant, which agrees with the constant
known from the energetic formulation. In Section~5, we present some
numerical results which confirm the a priori error estimates, as
given in Section~4. 
In Section~6, we finally draw some conclusions for future work.

\section{Energetic space-time boundary integral equation}
As in \cite{Aimi:2009}, we consider the Dirichlet boundary value problem
for the homogeneous wave equation in the one-dimensional spatial domain
$\Omega = (0,L)$ with zero initial conditions, and for a given
time horizon $T>0$,
\begin{equation}\label{DBVP}
  \left. \begin{array}{rcll}
           \partial_{tt} u(x,t) - \partial_{xx} u(x,t)
           & = & 0 & \mbox{for} \; (x,t) \in
                     Q := (0,L) \times (0,T), \\[1mm]
           u(x,0) = \partial_t u(x,t)_{|t=0}
           & = & 0 & \mbox{for} \; x \in (0,L),\\[1mm]
           u(0,t)
           & = & g_0(t) & \mbox{for} \; t \in (0,T), \\[1mm]
           u(L,t)
           & = & g_L(t) & \mbox{for} \; t \in (0,T).
         \end{array} \right \}
\end{equation}
In the one-dimensional case, the fundamental solution of the wave
equation is the Heaviside function
\[
U^*(x,t) = \frac{1}{2} \, H(t-|x|),
\]
and we can represent the solution $u$ of (\ref{DBVP})
by using the single layer potential
\[
  u(x,t) = (\widetilde{V}w)(x,t) = \frac{1}{2} \int_0^{t-|x|} w_0(s) \, ds +
  \frac{1}{2} \int_0^{t-|x-L|} w_L(s) \, ds \quad \mbox{for} \;
  (x,t) \in Q
\]
with the density functions $w=(w_0,w_L)$.
Note that for any function $z \colon \, (0,T) \to \R$, we set 
$z(t) = 0$ for $t<0$ or $t>T$ in the remainder of this work. 
To determine the yet unknown density functions $(w_0,w_L)$, we consider
the boundary integral equations for $x \to 0$,
\begin{equation}\label{BIE 0}
(V_0w)(t) :=
  \frac{1}{2} \int_0^t w_0(s) \, ds +
\frac{1}{2} \int_0^{t-L} w_L(s) \, ds = g_0(t) \quad \mbox{for} \;
t \in (0,T),
\end{equation}
and for $ x \to L$,
\begin{equation}\label{BIE L}
  (V_Lw)(t) :=
\frac{1}{2} \int_0^{t-L} w_0(s) \, ds +
  \frac{1}{2} \int_0^t w_L(s) \, ds = g_L(t) \quad \mbox{for} \; t \in (0,T).
\end{equation}
We write the boundary integral equations (\ref{BIE 0}) and (\ref{BIE L}) in
compact form, for $w=(w_0,w_L)$, as
\begin{equation}\label{BIE}
  (Vw)(t) =
  \begin{pmatrix}
           (V_0w)(t) \\
           (V_Lw)(t)
  \end{pmatrix}
  = \begin{pmatrix}
        V_{00} & V_{0L} \\ 
        V_{L0} & V_{LL}
    \end{pmatrix}
    \begin{pmatrix}
        w_0 \\ 
        w_L
    \end{pmatrix}
    (t)=
    \begin{pmatrix}
        g_0(t) \\ 
        g_L(t)
    \end{pmatrix}
    = g(t), \quad t \in (0,T) .
\end{equation}                         
In the energetic boundary element method \cite{Aimi:2009}, instead
of (\ref{BIE}), the time derivative of (\ref{BIE}) is considered,
\begin{equation}\label{BIE energetic}
\partial_t (Vw)(t) = \partial_t g(t) \quad \mbox{for} \; t \in (0,T).
\end{equation}
We introduce the related energetic bilinear form
\begin{eqnarray*}
  a(w,v) 
  &:=& \langle v , \partial_t V w \rangle_{L^2(\Sigma)} \\
  &=& \frac{1}{2} \int_0^T v_0(t) \frac{d}{dt} \int_0^t w_0(s) \, ds \, dt
         + \frac{1}{2} \int_0^T v_0(t) \frac{d}{dt}
         \int_0^{t-L} w_L(s) \, ds \, dt \\
  &&+ \frac{1}{2} \int_0^T v_L(t) \frac{d}{dt}
      \int_0^{t-L} w_0(s) \, ds \, dt
      + \frac{1}{2} \int_0^T v_L(t) \frac{d}{dt} \int_0^t w_L(s) \, ds \, dt \\
  &=& \frac{1}{2} \, \int_0^T v_0(t) \, w_0(t) \, dt +
        \frac{1}{2} \, \int_0^T v_0(t) \, w_L(t-L) \, dt \\
  &&+ \frac{1}{2} \, \int_0^T v_L(t) \, w_0(t-L) \, dt +
        \frac{1}{2} \, \int_0^T v_L(t) \, w_L(t) \, dt \, .
\end{eqnarray*}
When using both the Cauchy--Schwarz and H\"older inequality, we conclude
\begin{eqnarray*}
  |a(w,v)| & \leq & \frac{1}{2} \, \| v_0 \|_{L^2(0,T)} \| w_0 \|_{L^2(0,T)} +
                   \frac{1}{2} \, \| v_0 \|_{L^2(0,T)}
                   \| w_L \|_{L^2(0,T-L)} \\
          &&+ \frac{1}{2} \, \| v_L \|_{L^2(0,T)} \| w_0 \|_{L^2(0,T-L)}
              + \frac{1}{2} \, \| v_L \|_{L^2(0,T)} \| w_L \|_{L^2(0,T)} \\
          &\leq& \frac{1}{2} \, \| v_0 \|_{L^2(0,T)} \Big[
                   \| w_0 \|_{L^2(0,T)} + \| w_L \|_{L^2(0,T)} \Big] \\
          &&+ \frac{1}{2} \, \| v_L \|_{L^2(0,T)} \Big[
              \| w_0 \|_{L^2(0,T)} + \| w_L \|_{L^2(0,T)} \Big] \\
          &=& \frac{1}{2} \, \Big[
                \| v_0 \|_{L^2(0,T)} + \| v_L \|_{L^2(0,T)} \Big]
                \Big[ \| w_0 \|_{L^2(0,T)} + \| w_L \|_{L^2(0,T)} \Big] \\
          &\leq& \sqrt{
                \| v_0 \|_{L^2(0,T)}^2 + \| v_L \|_{L^2(0,T)}^2}
                   \sqrt{ \| w_0 \|^2_{L^2(0,T)} + \| w_L \|^2_{L^2(0,T)}} \\
          &=& \| v \|_{L^2(\Sigma)} \| w \|_{L^2(\Sigma)}
\end{eqnarray*}
for all $v=(v_0,v_L), w = (w_0,w_L) \in L^2(\Sigma) := L^2(0,T) \times L^2(0,T)$, where
\[
    \| z \|_{L^2(\Sigma)} := \left( \| z_0 \|_{L^2(0,T)}^2 + \| z_L \|_{L^2(0,T)}^2  \right)^{1/2} \quad \text{ for } z=(z_0,z_L) \in L^2(\Sigma).
\]
Moreover, the energetic bilinear form
$a(\cdot,\cdot)$ is also $L^2(\Sigma)$-elliptic, see
\cite[Theorem 2.1]{Aimi:2009}. For later reference, we will give a simplified
proof of this result. For this, we introduce
\begin{equation}\label{Def time slices}
n := \min \Big \{ m \in {\mathbb{N}} : T \leq m L \Big \},
\end{equation}
which is the number of time slices $T_j := ((j-1)L,j L)$ for $j =1,\ldots,n$
in the case $ T=nL$. In the case $ T < nL$, we define the last time slice
as $T_n:=((n-1)L,T)$, while all the others remain unchanged.

\begin{theorem}{\rm \cite[Theorem 2.1]{Aimi:2009}}\label{Theorem Aimi}
  For all $w \in L^2(\Sigma)$, there holds the ellipticity estimate
  \begin{equation}\label{Ellipticity V 0}
    a(w,w) = \langle w , \partial_t V w \rangle_{L^2(\Sigma)} \, \geq \,
     \sin^2 \frac{\pi}{2(n+1)} \, \| w \|_{L^2(\Sigma)}^2,
  \end{equation}
  where the number $n \in {\mathbb{N}}$ of time slices is defined in
  {\rm (\ref{Def time slices})}.
\end{theorem}

\proof{
  For $w = (w_0,w_L) \in L^2(\Sigma)$, we write
  \begin{align*}
    2 \, a(&w,w) \\
    &= \int_0^T [w_0(t)]^2 \, dt +
        \int_0^T w_0(t) \, w_L(t-L) \, dt 
  + \int_0^T w_L(t) \, w_0(t-L) \, dt +
        \int_0^T [w_L(t)]^2 \, dt \\
    &= \sum\limits_{j=1}^n \left[
        \| w_0 \|_{L^2(T_j)}^2
        + \int_{T_j} w_0(t) w_L(t-L) \, dt +
        \int_{T_j} w_L(t) w_0(t-L) \, dt 
        + \| w_L \|_{L^2(T_j)}^2 
        \right] .
  \end{align*}
  For $t \in T_1$, we have $t-L < 0$, and therefore $w_0(t-L)=w_L(t-L)=0$
  follows. For $j=2,\ldots,n-1$, we have, using the Cauchy--Schwarz inequality,
  \begin{multline*}
    \int_{T_j} w_0(t) w_L(t-L) \, dt
    \leq \left( \int_{T_j} [w_0(t)]^2 \, dt \right)^{1/2}
             \left( \int_{T_j} [w_L(t-L)]^2 \, dt \right)^{1/2} \\
    \leq \left( \int_{T_j} [w_0(t)]^2 \, dt \right)^{1/2}
             \left( \int_{T_{j-1}} [w_L(t)]^2 \, dt \right)^{1/2} =
             \| w_0 \|_{L^2(T_j)} \| w_L \|_{L^2(T_{j-1})} \, .
  \end{multline*}
  Correspondingly, for $j=n$ and $T_n = ((j-1)L,T)$, $T \leq nL$, we have
  \begin{align*}
    \int_{T_n} w_0(t) & w_L(t-L) \, dt
    \leq \left( \int_{T_n} [w_0(t)]^2 \, dt \right)^{1/2}
             \left( \int_{(n-1)L}^T [w_L(t-L)]^2 \, dt \right)^{1/2} \\
    &= \left( \int_{T_n} [w_0(t)]^2 \, dt \right)^{1/2}
        \left( \int_{(n-2)L}^{T-L} [w_L(t)]^2 \, dt \right)^{1/2} \\
    &\leq \left( \int_{T_n} [w_0(t)]^2 \, dt \right)^{1/2}
        \left( \int_{T_{n-1}} [w_L(t)]^2 \, dt \right)^{1/2}
        = \| w_0 \|_{L^2(T_n)} \| w_L \|_{L^2(T_{n-1})} \, .
  \end{align*}
  Hence, we conclude
  \begin{align*}
    2 \, a(w,w)
    \geq& \sum\limits_{j=1}^n \Big[
             \| w_0 \|^2_{L^2(T_j)} + \| w_L \|^2_{L^2(T_j)} \Big] \\
    &-
        \sum\limits_{j=2}^n \Big[
        \| w_0 \|_{L^2(T_j)} \| w_L \|_{L^2(T_{j-1})} +
        \| w_L \|_{L^2(T_j)} \| w_0 \|_{L^2(T_{j-1})} \Big] \\
    =& {\scriptsize \left(
            \begin{pmatrix}
                   1 & -\frac{1}{2} & & & & & \\
                   -\frac{1}{2} & 1 & - \frac{1}{2} & & & & \\
                     & - \frac{1}{2} & 1 & - \frac{1}{2} & & & \\
                     & & - \frac{1}{2} & 1 & - \frac{1}{2} & & \\
                   & & & \ddots & \ddots & \ddots & \\
                   & & & & - \frac{1}{2} & 1 & - \frac{1}{2} \\
                   & & & & & - \frac{1}{2} & 1
            \end{pmatrix}
            \begin{pmatrix}
                   \| w_0 \|_{L^2(T_1)} \\
                   \| w_L \|_{L^2(T_2)} \\
                   \| w_0 \|_{L^2(T_3)} \\
                   \| w_L \|_{L^2(T_4)} \\
                   \vdots  \\
                   \| w_0 \|_{L^2(T_{n-1})} \\
                   \| w_L \|_{L^2(T_n)} 
           \end{pmatrix},
           \begin{pmatrix}
                   \| w_0 \|_{L^2(T_1)} \\
                   \| w_L \|_{L^2(T_2)} \\
                   \| w_0 \|_{L^2(T_3)} \\
                   \| w_L \|_{L^2(T_4)} \\
                   \vdots  \\
                   \| w_0 \|_{L^2(T_{n-1})} \\
                   \| w_L \|_{L^2(T_n)} 
           \end{pmatrix} \right) } \\
    &+ {\scriptsize \left(
            \begin{pmatrix}
                   1 & -\frac{1}{2} & & & & & \\
                   -\frac{1}{2} & 1 & - \frac{1}{2} & & & & \\
                     & - \frac{1}{2} & 1 & - \frac{1}{2} & & & \\
                     & & - \frac{1}{2} & 1 & - \frac{1}{2} & & \\
                   & & & \ddots & \ddots & \ddots & \\
                   & & & & - \frac{1}{2} & 1 & - \frac{1}{2} \\
                   & & & & & - \frac{1}{2} & 1
            \end{pmatrix}
            \begin{pmatrix}
                   \| w_L \|_{L^2(T_1)} \\
                   \| w_0 \|_{L^2(T_2)} \\
                   \| w_L \|_{L^2(T_3)} \\
                   \| w_0 \|_{L^2(T_4)} \\
                   \vdots  \\
                   \| w_L \|_{L^2(T_{n-1})} \\
                   \| w_0 \|_{L^2(T_n)} 
            \end{pmatrix},
            \begin{pmatrix}
                   \| w_L \|_{L^2(T_1)} \\
                   \| w_0 \|_{L^2(T_2)} \\
                   \| w_L \|_{L^2(T_3)} \\
                   \| w_0 \|_{L^2(T_4)} \\
                   \vdots  \\
                   \| w_L \|_{L^2(T_{n-1})} \\
                   \| w_0 \|_{L^2(T_n)} 
           \end{pmatrix} \right) }
  \end{align*}
  and further,
  \begin{align*}
      a(w,w)
     &\geq \frac{\lambda_{\min}}{2} {\scriptsize \left[
        \left(         
        \begin{pmatrix}
                   \| w_0 \|_{L^2(T_1)} \\
                   \| w_L \|_{L^2(T_2)} \\
                   \| w_0 \|_{L^2(T_3)} \\
                   \| w_L \|_{L^2(T_4)} \\
                   \vdots  \\
                   \| w_0 \|_{L^2(T_{n-1})} \\
                   \| w_L \|_{L^2(T_n)} 
       \end{pmatrix},
       \begin{pmatrix}
                   \| w_0 \|_{L^2(T_1)} \\
                   \| w_L \|_{L^2(T_2)} \\
                   \| w_0 \|_{L^2(T_3)} \\
                   \| w_L \|_{L^2(T_4)} \\
                   \vdots  \\
                   \| w_0 \|_{L^2(T_{n-1})} \\
                   \| w_L \|_{L^2(T_n)} 
       \end{pmatrix} \right)
    +
    \left( \begin{pmatrix}
                   \| w_L \|_{L^2(T_1)} \\
                   \| w_0 \|_{L^2(T_2)} \\
                   \| w_L \|_{L^2(T_3)} \\
                   \| w_0 \|_{L^2(T_4)} \\
                   \vdots  \\
                   \| w_L \|_{L^2(T_{n-1})} \\
                   \| w_0 \|_{L^2(T_n)} 
            \end{pmatrix},
            \begin{pmatrix}
                   \| w_L \|_{L^2(T_1)} \\
                   \| w_0 \|_{L^2(T_2)} \\
                   \| w_L \|_{L^2(T_3)} \\
                   \| w_0 \|_{L^2(T_4)} \\
                   \vdots  \\
                   \| w_L \|_{L^2(T_{n-1})} \\
                   \| w_0 \|_{L^2(T_n)} 
           \end{pmatrix} \right) \right] } \\
    &= \frac{\lambda_{\min}}{2} \, \Big[ \| w_0 \|^2_{L^2(0,T)} +
          \| w_L \|^2_{L^2(0,T)} \Big],
  \end{align*}  
  where
  \[
\lambda_{\min} = 2 \, \sin^2 \frac{\pi}{2(n+1)}
  \]
  is the minimal eigenvalue of the involved matrix, which is related
  to the finite difference approximation of the Laplacian in one dimension.
}

\noindent
From the above properties, we conclude that
\[
  \partial_t V : L^2(\Sigma) \to L^2(\Sigma)
\]
defines an isomorphism.
Since the time derivative
\[
  \partial_t : H^1_{0,}(\Sigma) \to L^2(\Sigma)
\]
is also an isomorphism, e.g., \cite[Sect.~2.1]{SteinbachZank:2018}, so is
\begin{equation}\label{V mapping 1}
V : L^2(\Sigma) \to H^1_{0,}(\Sigma) \, .
\end{equation}
Note that, for $ u = (u_0,u_L) \in
H^1_{0,}(\Sigma) := H^1_{0,}(0,T) \times H^1_{0,}(0,T)$, we have
\[
  \| u \|^2_{H^1_{0,}(\Sigma)} := \| \partial_t u_0 \|^2_{L^2(0,T)} +
  \| \partial_t u_L \|^2_{L^2(0,T)} .
\]
For $\partial_t : H^1_{0,}(0,T) \to L^2(0,T)$, the inverse is given by
\[
  u(t) = (\partial_t^{-1} f)(t) = \int_0^t f(s) \, ds, \quad t \in (0,T),  
\]
with $f \in L^2(0,T)$, $u \in H^1_{0,}(0,T)$.
Analogously, for $\partial_t : H^1_{,0}(0,T) \to L^2(0,T)$, we find the
inverse as
\[
    u(t) = (\overline{\partial}_t^{-1} f)(t) = - \int_t^T f(s) \, ds, \quad t \in (0,T).
\]
For $w,v \in L^2(\Sigma)$ and $ u = Vw = (u_0,u_L) \in H^1_{0,}(\Sigma)$,
we therefore obtain
\[
  \langle \overline{\partial}_t^{-1} V w , v \rangle_{L^2(\Sigma)}
  =
  - \int_0^T \int_t^T u_0(s) \, ds \, v_0(t) \, dt
  - \int_0^T \int_t^T u_L(s) \, ds \, v_L(t) \, dt \, .
\]
For $ \ast \in \{ 0 ,L \}$ we compute
\begin{align*}
  - \int_0^T \int_t^T u_\ast(s) \, ds \, v_\ast(t) \, dt
  &= - \int_0^T \int_t^T u_\ast(s) \, ds \, \partial_t
        \int_0^t v_\ast(s) \, ds \, dt \\
  =&
        \left. - \int_t^T u_\ast(s) \, ds \int_0^t v_\ast(s) \, ds \right|_0^T
      + \int_0^T \partial_t \int_t^T u_\ast(s) \, ds
      \int_0^t v_\ast(s) \, ds \, dt \\
  =& - \int_0^T u_\ast(t) \int_0^t v_\ast(s) \, ds \, ,
\end{align*}
i.e.,
\[
  \langle \overline{\partial}_t^{-1} Vw , v \rangle_{L^2(\Sigma)} =
  - \langle V w , \partial_t^{-1} v \rangle_{L^2(\Sigma)} \, .
\]
On the other hand, for $ z_0 = \partial_t^{-1} w_0$ we have
$ w_0 = \partial_t z_0$, and hence
\[
  \int_0^t w_0(s) \, ds =
  \int_0^t \partial_s z_0(s) \, ds = z_0(t) = \partial_t \int_0^t z_0(s) \, ds . 
\]
With this, we conclude
\[
  \langle \overline{\partial}_t^{-1} Vw , v \rangle_{L^2(\Sigma)} =
  - \langle V \partial_t \partial_t^{-1} w ,
  \partial_t^{-1} v \rangle_{L^2(\Sigma)} =
  - \langle \partial_t V \partial_t^{-1} w , \partial_t^{-1} v \rangle_{L^2(\Sigma)} = -a( \partial_t^{-1} w , \partial_t^{-1} v),
\]
and, in particular for $v=w$, Theorem~\ref{Theorem Aimi} gives
\[
- \langle \overline{\partial}_t^{-1} V w , w \rangle_{L^2(\Sigma)} =
  \langle \partial_t V \partial_t^{-1} w , \partial_t^{-1}w \rangle_{L^2(\Sigma)}
  \nonumber
  \geq \sin^2 \frac{\pi}{2(n+1)} \,
           \| \partial_t^{-1} w \|_{L^2(\Sigma)}^2 .
\]
For $\ast \in \{ 0 , L \}$, we define
\[
z_\ast(t) = (\partial_t^{-1} w_\ast)(t) = \int_0^t w_\ast(s) \, ds, \quad t \in (0,T), 
\]
to compute
\begin{align*}
  \| \partial_t^{-1} w_\ast \|^2_{L^2(0,T)}
  &= \| z_\ast \|^2_{L^2(0,T)} =
        \int_0^T z_\ast(t) \, z_\ast(t) \, dt =
        - \int_0^T \partial_t \int_t^T z_\ast(s) \, ds \, z_*(t) \, dt \\
  &= \left. - \int_t^T z_\ast(s) \, ds \, z_\ast(t) \right|_0^T +
        \int_0^T \int_t^T z_*(s) \, ds \; \partial_t z_\ast(t) \, dt \\
  &= \int_0^T v_*(t) \, w_*(t) \, dt,
\end{align*}
where
\[
  v_*(t) = \int_t^T z_*(s) \, ds \quad \text{ for } t \in (0,T), \quad
  \partial_t v_* = -z_*, \quad v_* \in H^1_{,0}(0,T) .
\]
From this, we conclude
\[
  \| \partial_t^{-1} w_* \|_{L^2(0,T)} =
  \frac{| \langle w_* , v_* \rangle_{(0,T)}| }{\| \partial_t v_* \|_{L^2(0,T)}}
  \leq \sup\limits_{0 \neq \phi \in H^1_{,0}(0,T)}
  \frac{|\langle w_* , \phi \rangle_{(0,T)}|}{\| \partial_t \phi \|_{L^2(0,T)}}
  = \| w_* \|_{[H^1_{,0}(0,T)]'} .
\]
Indeed, we have
\[
  \| \partial_t^{-1} w_* \|_{L^2(0,T)} 
  = \| w_* \|_{[H^1_{,0}(0,T)]'} ,
\]
and therefore,
\begin{equation}
  - \langle \overline{\partial}_t^{-1} V w , w \rangle_{L^2(\Sigma)}
  \geq \sin^2 \frac{\pi}{2(n+1)} \, \| w \|_{[H^1_{,0}(\Sigma)]'}^2 .
        \label{Ellipticity V -1}
\end{equation}
In fact, by the density of $L^2(\Sigma)$ in $[H^1_{,0}(\Sigma)]'$, the operator
\[
- \overline{\partial}_t^{-1} V : [H^1_{,0}(\Sigma)]' \to H^1_{,0}(\Sigma) 
\]
defines an isomorphism, and so does
\begin{equation}\label{V mapping 2}
V : [H^1_{,0}(\Sigma)]' \to L^2(\Sigma) .
\end{equation}
For the single layer boundary integral operator $V$, we have obtained
the mapping properties (\ref{V mapping 1}) and (\ref{V mapping 2}),
respectively. When applying an interpolation argument, this gives that
\begin{equation*}
    V : [H^{1/2}_{,0}(\Sigma)]' \to H^{1/2}_{0,}(\Sigma) 
\end{equation*}
is an isomorphism as well, where the Sobolev space $H^{1/2}_{0,}(\Sigma) = 
H^{1/2}_{0,}(0,T) \times H^{1/2}_{0,}(0,T)$ is endowed with the Hilbertian norm
\[
 \| z \|_{H^{1/2}_{0,}(\Sigma)} := \left( \| z_0 \|_{H^{1/2}_{0,}(0,T)}^2 
 + \| z_L \|_{H^{1/2}_{0,}(0,T)}^2\right)^{1/2} 
    \quad \text{ for } z=(z_0,z_L) \in H^{1/2}_{0,}(\Sigma)
\]
and analogously, the Sobolev space $H^{1/2}_{,0}(\Sigma)$ is introduced.
Hence, we conclude the inf-sup stability condition
\begin{equation}\label{inf sup V}
  c_S \, \| w \|_{[H^{1/2}_{,0}(\Sigma)]'} \leq
  \sup\limits_{0 \neq v \in [H^{1/2}_{0,}(\Sigma)]'}
  \frac{|\langle V w , v \rangle_\Sigma|}{\| v \|_{[H^{1/2}_{0,}(\Sigma)]'}}
  \quad \mbox{for all} \; w \in [H^{1/2}_{,0}(\Sigma)]'
\end{equation}
with a constant $c_S > 0.$
In fact, \eqref{inf sup V} corresponds to the
inf-sup condition in \cite[Theorem 5.7]{SteinbachUrzua:2021}, where
the test space is slightly larger than used in (\ref{inf sup V}). But
we will show that
$V : [H^{1/2}_{,0}(\Sigma)]' \to H^{1/2}_{0,}(\Sigma)$ in combination with
a modified Hilbert transformation \cite{SteinbachZank:2018, SteinbachZank:2021, Zank2021Exact} even satisfies
an ellipticity estimate similar as in (\ref{Ellipticity V 0}).

\section{A modified Hilbert transformation}
For $ u \in L^2(0,T)$, we consider the Fourier series
\[
  u(t) = \sum_{k=0}^\infty u_k \, \sin \left( \left( \frac{\pi}{2} + k\pi
    \right) \frac{t}{T} \right), \quad
    u_k = \frac{2}{T} \int_0^T u(t) \,
    \sin \left( \left( \frac{\pi}{2} + k\pi
    \right) \frac{t}{T} \right) \, dt ,
\]
\[
  u(t) = \sum_{k=0}^\infty \overline{u}_k \,
  \cos \left( \left( \frac{\pi}{2} + k\pi
    \right) \frac{t}{T} \right), \quad
    \overline{u}_k = \frac{2}{T} \int_0^T u(t) \,
    \cos \left( \left( \frac{\pi}{2} + k\pi
    \right) \frac{t}{T} \right) \, dt .
\]
From \cite[Lemma 2.1]{SteinbachZank:2018}, we have
\begin{equation*}
  \| u \|^2_{[H^{1/2}_{,0}(0,T)]'} = \frac{T^2}{2} \,
  \sum\limits_{k=0}^\infty \left( \frac{\pi}{2} + k\pi \right)^{-1}
  \overline{u}_k^2 \, .
\end{equation*}
As in \cite{SteinbachZank:2018}, we introduce the transformation
operator ${\mathcal{H}}_T : L^2(0,T) \to L^2(0,T)$ as
\begin{equation}\label{Def H}
  {\mathcal{H}}_T u(t) := \sum\limits_{k=0}^\infty u_k \, 
  \cos \left( \left( \frac{\pi}{2} + k\pi
    \right) \frac{t}{T} \right), \quad t \in (0,T),
\end{equation}
which is norm preserving and bijective. By construction, we have that the transformation operator ${\mathcal{H}}_T : H^{1/2}_{0,}(0,T) \to H^{1/2}_{,0}(0,T)$ is also an isometric isomorphism, and
\[
  \langle \partial_t u , {\mathcal{H}}_T u \rangle_{(0,T)} =
  \| u \|^2_{H^{1/2}_{0,}(0,T)} \quad \mbox{for all} \;
  u \in H^{1/2}_{0,}(0,T) .
\]
Note that $H^{1/2}_{0,}(0,T) := [H^1_{0,}(0,T),L^2(0,T)]_{1/2}$ is constructed
by interpolation, where $H^1_{0,}(0,T) := \{ v \in H^1(0,T):v(0)=0 \}$.
In the same way, we define $H^{1/2}_{,0}(0,T)$ but with zero condition at the
final time $t=T$. It is easy to see that
\begin{equation}\label{Cauchy Schwarz Hilbert}
  |\langle \partial_t u , {\mathcal{H}}_T z \rangle_{(0,T)} | \leq
  \| u \|_{H^{1/2}_{0,}(0,T)} \| z \|_{H^{1/2}_{0,}(0,T)} \quad
  \mbox{for all} \; u,z \in H^{1/2}_{0,}(0,T).
\end{equation}
The transformation operator ${\mathcal{H}}_T$, as defined in
(\ref{Def H}), allows a closed representation, see
\cite[Lemma 2.8]{SteinbachZank:2018}, which generalizes the well-known
Hilbert transformation, e.g., \cite{ButzerTrebels:1968}.
Moreover, following \cite[Eqn. (2.5)]{SteinbachZank:2021}
we conclude the following representation for $u,z \in H^{1}_{0,}(0,T)$,
\begin{equation*}
  \langle \partial_t u , {\mathcal{H}}_T z \rangle_{(0,T)}
  = - \frac{1}{\pi} \int_0^T \partial_t u(t) \int_0^T \ln \left[
    \tan \frac{\pi(s+t)}{4T} \tan \frac{\pi|t-s|}{4T} \right]
  \, \partial_s z(s) \, ds \, dt \, .
\end{equation*}
This representation also allows for an efficient evaluation of the
bilinear form $\langle \partial_t u , {\mathcal{H}}_T z \rangle_{(0,T)}$
by using hierarchical matrices, see \cite{SteinbachZank:2021} for a
more detailed discussion.

\section{A space-time approach in energy spaces}
Instead of the boundary integral equation \eqref{BIE energetic},
we may replace the application of the time derivative by the
modified Hilbert transformation
${\mathcal{H}}_T : H^{1/2}_{0,}(\Sigma) \to H^{1/2}_{,0}(\Sigma)$,
i.e., we consider the boundary integral equation to find
$ w \in [H^{1/2}_{,0}(\Sigma)]'$ such that
\begin{equation*}
  {\mathcal{H}}_T V w = {\mathcal{H}}_T g \quad \mbox{in} \;
  [H^{1/2}_{,0}(\Sigma)]',
\end{equation*}
where $g \in H^{1/2}_{0,}(\Sigma)$ is a given Dirichlet datum. 
The related bilinear form is given as
\[
  a_{{\mathcal{H}}_T}(w,v) :=
  \langle v , {\mathcal{H}}_T V w \rangle_{\Sigma} \quad
  \mbox{for all} \; v,w \in [H^{1/2}_{,0}(\Sigma)]' .
\]
Recall that for $ u = (u_0,u_L) \in H^{1/2}_{0,}(\Sigma)$, we have
\[
  \partial_t u = ( \partial_t u_0 , \partial_t u_L) =
  (v_0,v_L) =: v \in [H^{1/2}_{,0}(\Sigma)]' \, ,
\]
satisfying
\[
\| u \|_{H^{1/2}_{0,}(\Sigma)} = \| v \|_{[H^{1/2}_{,0}(\Sigma)]'} \, .
\]
For $v = \partial_t u$, $w = \partial_t z$ with $u, z \in H^{1}_{0,}(\Sigma)$, we can write
\begin{eqnarray*}
  a_{{\mathcal{H}}_T}(w,v)
  & = & \frac{1}{2} \int_0^T v_0(t) \, {\mathcal{H}}_T \left(
        \int_0^t w_0(s) \, ds + \int_0^{t-L} w_L(s) \, ds \right) dt \\
  & & + \frac{1}{2} \int_0^T v_L(t) \, {\mathcal{H}}_T \left(
      \int_0^{t-L} w_0(s) \, ds + \int_0^t w_L(s) \, ds \right) dt \\
  & = & \frac{1}{2} \, 
        \Big[ \langle \partial_t u_0 , {\mathcal{H}}_T (z_0+z_L(\cdot-L))
        \rangle_{(0,T)} +
        \langle \partial_t u_L , {\mathcal{H}}_T (z_0(\cdot - L)+z_L)
        \rangle_{(0,T)} \Big] \, .
\end{eqnarray*}
When using \eqref{Cauchy Schwarz Hilbert}, we obtain
\begin{align*}
  |a_{{\mathcal{H}}_T}(&w,v)| \\
  &\leq \frac{1}{2} \, \Big[
           \| u_0 \|_{H^{1/2}_{0,}(0,T)}
           \| z_0 + z_L(\cdot - L) \|_{H^{1/2}_{0,}(0,T)} +
           \| u_L \|_{H^{1/2}_{0,}(0,T)}
      \| z_0(\cdot - L) + z_L \|_{H^{1/2}_{0,}(0,T)} \Big] \\
  &\leq \frac{1}{2} \, \Big[
      \| u_0 \|_{H^{1/2}_{0,}(0,T)} + \| u_L \|_{H^{1/2}_{0,}(0,T)} \Big]
      \Big[
      \| z_0 \|_{H^{1/2}_{0,}(0,T)} + \| z_L \|_{H^{1/2}_{0,}(0,T)} \Big] \\
  &\leq \sqrt{
      \| u_0 \|_{H^{1/2}_{0,}(0,T)}^2 + \| u_L \|^2_{H^{1/2}_{0,}(0,T)}}
      \sqrt{
      \| z_0 \|^2_{H^{1/2}_{0,}(0,T)} + \| z_L \|^2_{H^{1/2}_{0,}(0,T)}}
  \\
  &= \| u \|_{H^{1/2}_{0,}(\Sigma)} \| z \|_{H^{1/2}_{0,}(\Sigma)} \\
  &= \| v \|_{[H^{1/2}_{,0}(\Sigma)]'} \| w \|_{[H^{1/2}_{,0}(\Sigma)]'} 
\end{align*}
for all $v,w \in L^2(\Sigma)$, i.e., the density of $L^2(\Sigma)$ in
$[H^{1/2}_{,0}(\Sigma)]'$ yields the boundedness of the
bilinear form $a_{{\mathcal{H}}_T}(\cdot,\cdot)$.

\begin{lemma}\label{Lemma HV}
For $ w \in [H^{1/2}_{,0}(\Sigma)]'$, there holds
\begin{equation}\label{HV elliptic}
  a_{{\mathcal{H}}_T}(w,w) = \langle {\mathcal{H}}_T V w , w \rangle_\Sigma
  \geq 
  \frac{1}{2} \left( 1 - \frac{1}{2} \,
        \sup\limits_{m \in {\mathbb{N}}} \sqrt{\lambda_{\max}(C_m)} \right)
        \, \| w \|^2_{[H^{1/2}_{,0}(\Sigma)]'} ,
\end{equation}
where $\lambda_{\max}(C_m)$ is the maximal eigenvalue of a symmetric matrix
$C_m \in \R^{(m+1)\times (m+1)}$. In the case of $T \leq L$, the matrix $C_m$ is
the zero matrix, i.e., $\lambda_{\max}(C_m) = 0$. However, in the case $T > L$, the matrix $C_m$ is defined by the entries
\[
  c_{\ell i} = \sum\limits_{k=0}^\infty b_{k\ell} b_{ki} \quad \mbox{for} \;
  \ell,i=0,\ldots,m,
\]
\[
b_{kk} = 2 \, \left( 1-\frac{L}{T} \right)
        \cos \left(
        \left( \frac{\pi}{2} + k\pi \right) \frac{L}{T}
        \right) \quad \mbox{for} \; k \in {\mathbb{N}}_0,
\]
\[
b_{k\ell} = \frac{4}{\pi} \frac{\sqrt{2k+1}\sqrt{2\ell+1}}
     {(k+\ell+1)(k-\ell)}
        \cos \left( (k+\ell+1) \frac{\pi}{2} \frac{L}{T} \right)
        \sin \left( (\ell-k) \frac{\pi}{2} \frac{L}{T} \right)
\]
for $k,\ell \in {\mathbb{N}}_0$, $k-\ell=2j \neq 0$, $j \in {\mathbb{Z}}$,
and $b_{k\ell} = 0$ else.
\end{lemma}

\proof
{
For $w=(w_0,w_L) \in L^2(\Sigma)$, we consider the
  Fourier series
\[
  w_0(t) = \sum\limits_{k=0}^\infty \overline{w}_{0,k}
  \cos \left( \left( \frac{\pi}{2} + k\pi \right) \frac{t}{T}\right), \quad
  \overline{w}_{0,k} = \frac{2}{T} \int_0^T w_0(t) \,
  \cos \left( \left( \frac{\pi}{2} + k\pi \right) \frac{t}{T} \right) dt,
\]
\[
  w_L(t) = \sum\limits_{k=0}^\infty \overline{w}_{L,k}
  \cos \left( \left( \frac{\pi}{2} + k\pi \right) \frac{t}{T}\right), \quad
  \overline{w}_{L,k} = \frac{2}{T} \int_0^T w_L(t) \,
  \cos \left( \left( \frac{\pi}{2} + k\pi \right) \frac{t}{T} \right) dt.
\]
In the case $T \leq L$, we explicitly compute
\begin{equation*}
 \langle {\mathcal{H}}_T V w , w \rangle_{L^2(\Sigma)} 
 =\frac{T^2}{2} \sum\limits_{k=0}^\infty \frac{\overline{w}_{0,k}^2 + \overline{w}_{L,k}^2}{(2k+1)\pi}
 = \frac 1 2 \Big( \| w_0 \|^2_{[H^{1/2}_{,0}(0,T)]'} + \| w_L \|^2_{[H^{1/2}_{,0}(0,T)]'} \Big),
\end{equation*}
since there are no coupling terms.

In the case $T > L$, we have the representation
\begin{align*}
  \langle {\mathcal{H}}_T V& w , w \rangle_{L^2(\Sigma)}
   \\
  =& \,
  \frac{T^2}{2} \sum\limits_{k=0}^\infty
  \frac{\overline{w}_{0,k}^2 + \overline{w}_{L,k}^2}{(2k+1)\pi}
  + \frac{T^2}{2}
  \sum\limits_{k=0}^\infty 
        \overline{w}_{0,k} \overline{w}_{L,k}
        \frac{2}{(2k+1)\pi} 
        \left( 1-\frac{L}{T} \right)
        \cos \left(
        \left( \frac{\pi}{2} + k\pi \right) \frac{L}{T}
        \right) \\
  &+ \frac{T^2}{2} \sum\limits_{k-\ell=2j \neq 0} \overline{w}_{0,\ell} \overline{w}_{L,k}
        \frac{4}{\pi^2} \frac{1}{(k+\ell+1)(k-\ell)}
        \cos \left( (k+\ell+1) \frac{\pi}{2} \frac{L}{T} \right)
     \sin \left( (\ell-k) \frac{\pi}{2} \frac{L}{T} \right) \\
  =& \,
  \frac{T^2}{2} \sum\limits_{k=0}^\infty
  \Big[ \widehat{w}_{0,k}^2 + \widehat{w}_{L,k}^2 \Big]
  + \frac{T^2}{2}
  \sum\limits_{k=0}^\infty 
        2 \widehat{w}_{0,k} \widehat{w}_{L,k}
        \left( 1-\frac{L}{T} \right)
        \cos \left(
        \left( \frac{\pi}{2} + k\pi \right) \frac{L}{T}
        \right) \\
  &+ \frac{T^2}{2} \sum\limits_{k-\ell=2j \neq 0} \widehat{w}_{0,\ell}
     \widehat{w}_{L,k}
     \frac{4}{\pi} \frac{\sqrt{2k+1}\sqrt{2\ell+1}}
     {(k+\ell+1)(k-\ell)}
        \cos \left( (k+\ell+1) \frac{\pi}{2} \frac{L}{T} \right)
        \sin \left( (\ell-k) \frac{\pi}{2} \frac{L}{T} \right),
\end{align*}
where
\[
  \widehat{w}_{0,k} = \frac{\overline{w}_{0,k}}{\sqrt{(2k+1)\pi}}, \quad
  \widehat{w}_{L,k} = \frac{\overline{w}_{L,k}}{\sqrt{(2k+1)\pi}} .
\]
When using the coefficients $b_{k\ell}$, we write the above result as
\[
  \langle {\mathcal{H}}_T V w , w \rangle_{L^2(\Sigma)}
  = \frac{T^2}{2} \left(
    \sum\limits_{k=0}^\infty \Big[ \widehat{w}^2_{0,\ell} +
    \widehat{w}^2_{L,k} \Big] + \sum\limits_{k=0}^\infty
    \sum\limits_{\ell=0}^\infty b_{k\ell} \widehat{w}_{0,\ell} \widehat{w}_{L,k}
        \right) \, .
\]
Following \cite[Chapter VIII]{Hardy:1952}, we consider the forms
\[
  B(\widehat{w}_0, \widehat{w}_L) :=
  \sum\limits_{k=0}^\infty
  \sum\limits_{\ell=0}^\infty b_{k\ell} \widehat{w}_{0,\ell} \widehat{w}_{L,k},
  \quad
  B_m(\widehat{w}_0, \widehat{w}_L) :=
  \sum\limits_{k=0}^m
  \sum\limits_{\ell=0}^m b_{k\ell} \widehat{w}_{0,\ell} \widehat{w}_{L,k} \, ,
\]
and for the latter we estimate
\begin{align*}
  \Big|  B_m(\widehat{w}_0, \widehat{w}_L) \Big| = \left| \sum\limits_{k=0}^m
  \sum\limits_{\ell=0}^m b_{k\ell} \widehat{w}_{0,\ell} \widehat{w}_{L,k}
        \right| 
  &\leq \left[ \sum\limits_{k=0}^m \widehat{w}_{L,k}^2 \right]^{1/2}
           \left[ \sum\limits_{k=0}^m \left(
           \sum\limits_{\ell=0}^m b_{k\ell} \widehat{w}_{0,\ell} \right)^2
        \right]^{1/2} \\
  & \leq \left[ \sum\limits_{k=0}^m \widehat{w}_{L,k}^2 \right]^{1/2}
           \left[ \sum\limits_{k=0}^\infty \left(
           \sum\limits_{\ell=0}^m b_{k\ell} \widehat{w}_{0,\ell} \right)^2
           \right]^{1/2} \, .
\end{align*}
Hence, it remains to consider
\begin{align*}
  \sum\limits_{k=0}^\infty \left(
  \sum\limits_{\ell=0}^m b_{k\ell} \widehat{w}_{0,\ell} \right)^2
  &= \sum\limits_{\ell=0}^m \sum\limits_{j=0}^m
        \left( \sum\limits_{k=0}^\infty b_{k\ell} b_{kj} \right)
        \widehat{w}_{0,\ell} \widehat{w}_{0,j} \\[1mm]
  &= \sum\limits_{\ell=0}^m \sum\limits_{j=0}^m c_{\ell j}
        \widehat{w}_{0,\ell} \widehat{w}_{0,j} 
        \, \leq \, \lambda_{\max}(C_m) \sum\limits_{\ell=0}^m
        \widehat{w}_{0,\ell}^2 \, .
\end{align*}
From this, we conclude
\begin{align*}
  \Big|  B_m(\widehat{w}_0, \widehat{w}_L) \Big|
  &\leq \sqrt{\lambda_{\max}(C_m)} 
           \left[ \sum\limits_{k=0}^m \widehat{w}_{L,k}^2 \right]^{1/2}
           \left[ 
           \sum\limits_{\ell=0}^m \widehat{w}_{0,\ell}^2 
           \right]^{1/2} \\
  &\leq \sup\limits_{m \in {\mathbb{N}}} \sqrt{\lambda_{\max}(C_m)} 
           \left[ \sum\limits_{k=0}^\infty \widehat{w}_{L,k}^2 \right]^{1/2}
           \left[ 
           \sum\limits_{\ell=0}^\infty \widehat{w}_{0,\ell}^2 
           \right]^{1/2} \\
  &\leq \frac{1}{2} \,
           \sup\limits_{m \in {\mathbb{N}}} \sqrt{\lambda_{\max}(C_m)} 
           \left( \sum\limits_{k=0}^\infty \widehat{w}_{L,k}^2 +
           \sum\limits_{\ell=0}^\infty \widehat{w}_{0,\ell}^2 
           \right)
\end{align*}
for all $m \in {\mathbb{N}}$, and therefore
\[
  \Big|  B(\widehat{w}_0, \widehat{w}_L) \Big| \leq
  \frac{1}{2} \,
           \sup\limits_{m \in {\mathbb{N}}} \sqrt{\lambda_{\max}(C_m)} 
           \left( \sum\limits_{k=0}^\infty \widehat{w}_{L,k}^2 +
           \sum\limits_{\ell=0}^\infty \widehat{w}_{0,\ell}^2 
           \right)
\]
follows. With this, we finally obtain
\begin{align*}
  \langle {\mathcal{H}}_T V w , w \rangle_{L^2(\Sigma)}
  &\geq
  \frac{T^2}{2} \left( 1 - \frac{1}{2} \,
    \sup\limits_{m \in {\mathbb{N}}} \sqrt{\lambda_{\max}(C_m)} \right)
  \sum\limits_{k=0}^\infty \Big[ \widehat{w}_{0,k}^2 + \widehat{w}_{L,k}^2
           \Big]\\
  &=
  \frac{T^2}{2} \left( 1 - \frac{1}{2} \,
    \sup\limits_{m \in {\mathbb{N}}} \sqrt{\lambda_{\max}(C_m)} \right)
  \sum\limits_{k=0}^\infty \frac{\overline{w}_{0,k}^2 + \overline{w}_{L,k}^2}{(2k+1) \pi} \\
  &=
  \frac{T^2}{4} \left( 1 - \frac{1}{2} \,
    \sup\limits_{m \in {\mathbb{N}}} \sqrt{\lambda_{\max}(C_m)} \right)
        \sum\limits_{k=0}^\infty \frac{\overline{w}_{0,k}^2 + \overline{w}_{L,k}^2}
        {\frac{\pi}{2}+k\pi} \\
  &= \frac{1}{2} \left( 1 - \frac{1}{2} \,
        \sup\limits_{m \in {\mathbb{N}}} \sqrt{\lambda_{\max}(C_m)} \right)
        \Big( \| w_0 \|^2_{[H^{1/2}_{,0}(0,T)]'} +
        \| w_L \|^2_{[H^{1/2}_{,0}(0,T)]'} \Big),
\end{align*}
as stated.
In both cases $T\leq L$ or $T>L$, the density of $L^2(\Sigma)$ in
$[H^{1/2}_{,0}(\Sigma)]'$ yields the assertion.
}

\begin{proposition}\label{Proposition C}
  Numerical results indicate that
  \[
    \sup\limits_{m \in {\mathbb{N}}} \sqrt{\lambda_{\max}(C_m)} =
    2 - 4 \sin^2 \left( \frac{\pi}{2(n+1)} \right),
  \]
  where $n$ is given in {\rm (\ref{Def time slices})}. Indeed, for
  $L=1$, $T \in [1,20]$ and $m =20000$, the related results
  are given in Figure~\ref{Fig:Results}. Then, the ellipticity estimate
  {\rm (\ref{HV elliptic})} becomes
  \begin{equation}\label{HV elliptic final}
    a_{{\mathcal{H}}_T}(w,w) = \langle {\mathcal{H}}_T V w , w \rangle_\Sigma
  \geq  \sin^2 \left( \frac{\pi}{2(n+1)} \right)
  \, \| w \|^2_{[H^{1/2}_{,0}(\Sigma)]'} \quad \mbox{for all} \;
  w \in [H^{1/2}_{,0}(\Sigma)]' ,
  \end{equation}
  where the ellipticity constant is the same as in
  {\rm (\ref{Ellipticity V 0})}, and in {\rm (\ref{Ellipticity V -1})},
  respectively. Hence, we can think of {\rm (\ref{HV elliptic final})} being
  an interpolation of the ellipticity estimates
  {\rm (\ref{Ellipticity V 0})} and {\rm (\ref{Ellipticity V -1})}.

\begin{figure}[ht]
\begin{center}
     \begin{tikzpicture}[font=\scriptsize,scale=1.4]
         \begin{axis}
             [
                 grid=major, grid style=dashed,
                 xtick={1,2,4,...,20},
                 xmin=0.5, xmax=20.5,
                 ymin=-0.1, ymax=2.1,
                 legend entries={$\sqrt{\lambda_{\max}(C_m)}$,
                   $2 - 4 \sin^2 \Big( \frac{\pi}{2(n+1)} \Big)$},
                 legend pos = {south east}, mark size=1pt
             ],
             \addplot[only marks, color=red, mark=*, mark
options={scale=0.7}] plot file {KonstanteAimi.tex};
             \addplot[only marks, color=blue, mark=+, mark
options={scale=1.2}] plot file {Eigenwerte.tex};
         \end{axis}
     \end{tikzpicture}
     \caption{Numerical evaluation of $\sqrt{\lambda_{\max}(C_m)}$ for
       $L=1$, $T \in [1,20]$, $m = 20000$.}\label{Fig:Results}
\end{center}
\end{figure}
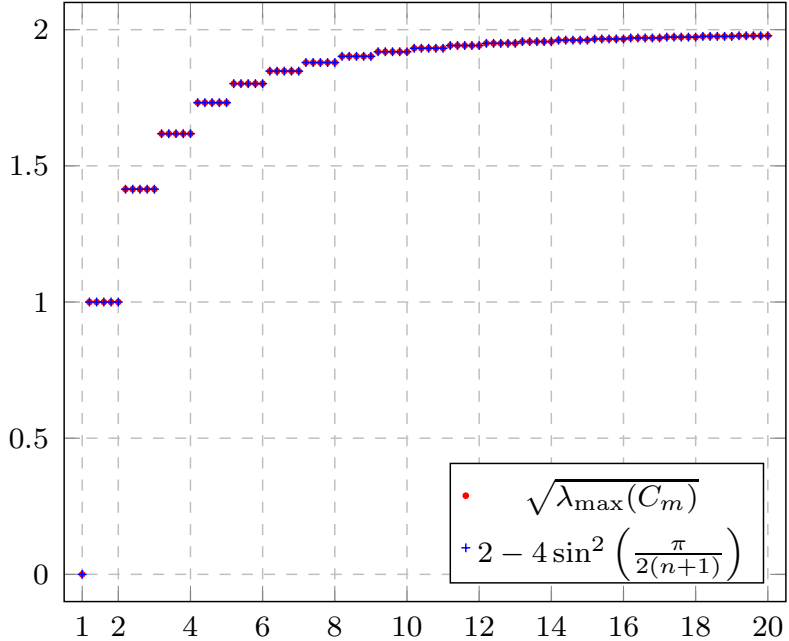

\end{proposition}

With the above results, we conclude unique solvability of the variational
formulation to find $w \in [H^{1/2}_{,0}(\Sigma)]'$ such that
\begin{equation}\label{VF}
  \langle v , {\mathcal{H}}_T V w \rangle_\Sigma =
  \langle v , {\mathcal{H}}_T g \rangle_\Sigma \quad
  \mbox{for all} \; v \in [H^{1/2}_{,0}(\Sigma)]',
\end{equation}
where $g \in H^{1/2}_{0,}(\Sigma)$ is a given Dirichlet datum. Let
$W_h \subset [H^{1/2}_{,0}(\Sigma)]'$ be some boundary element space,
e.g., of piecewise constant basis functions, which are defined with
respect to some decomposition of the lateral boundaries
$ \{ 0 \} \times (0,T)$ and $ \{ L \} \times (0,T)$, respectively.
The space-time Galerkin boundary element formulation of (\ref{VF}) is:
Find $w_h \in W_h$ such that
\begin{equation*}
  \langle v_h , {\mathcal{H}}_T V w_h \rangle_\Sigma =
  \langle v_h , {\mathcal{H}}_T g \rangle_\Sigma \quad
  \mbox{for all} \; v_h \in W_h .
\end{equation*}
When assuming $ w \in H^s(\Sigma)$ for some $ s \in [0,1]$ and
using standard arguments, e.g., \cite{Steinbach:2008}, we derive
an a priori error estimate in the energy norm,
\begin{equation*}
  \| w - w_h \|_{[H^{1/2}_{,0}(\Sigma)]'} \leq c \, h^{s+\frac{1}{2}} \,
  \|w\|_{H^s(\Sigma)} .
\end{equation*}
Moreover, using an inverse inequality, we also obtain an error estimate
in $L^2(\Sigma)$,
\begin{equation}\label{error L2}
  \| w - w_h \|_{L^2(\Sigma)} \leq c \, h^s \, \|w\|_{H^s(\Sigma)} .
\end{equation}

\section{Numerical results}
Instead of the boundary integral equation (\ref{BIE}) of the indirect
approach, we consider, as in \cite{SteinbachZank:2016}, the
boundary integral equation of the direct approach
\begin{equation} \label{Num:DirektBIE}
    V w = (\frac{1}{2}I+K)g \quad \text{ on } \Sigma,
\end{equation}
including the double layer boundary integral operator $K$ on the right hand 
side. In this case, the unknown $w$ is the spatial normal derivative 
$\partial_{n_x} u$ of the solution $u$ of \eqref{DBVP}.

For a boundary element approximation, consider a decomposition of the lateral boundary
\begin{equation*}
	\Sigma = \bigcup_{i=1}^{N_0+N_L} \overline{\tau}_i
\end{equation*}
into $N_0+N_L$ boundary elements $\tau_i$ with maximal mesh size $h = \max_{i} |\tau_i|$. Here, $N_0$ is the number of boundary elements for the boundary $\{ 0 \} \times (0,T)$ and  $N_L$ is the number of boundary elements for the boundary $\{ L \} \times (0,T)$. The conforming ansatz space of piecewise constant functions
\begin{equation*}
		S_{h}^0 (\Sigma) := S_{h_0}^0 (0,T) \times S_{h_L}^0 (0,T) \subset [H^{1/2}_{,0}(\Sigma)]'
\end{equation*}
is used to define an approximate solution $w_h \in S_h^0(\Sigma).$ Then, the Galerkin discretization of \eqref{Num:DirektBIE} to find $w_h \in S_h^0(\Sigma)$ such that
\begin{equation} \label{Num:direktBEM}
 \langle v_h, \mathcal H_T V w_h \rangle_{L^2(\Sigma)}  = \langle v_h, \mathcal H_T (\frac{1}{2}I+K) Q_h g \rangle_{L^2(\Sigma)} \quad \text{ for all } v_h \in S_h^0(\Sigma) 
\end{equation}
is equivalent to the global linear system
\begin{equation} \label{Num:LGS}
	V_h \underline w = \underline g
\end{equation}
with the related system matrix $V_h \in \R^{(N_0 + N_L) \times (N_0 + N_L)},$ 
the right-hand side $\underline g \in \R^{N_0 + N_L}$ and the vector of unknown 
coefficients $\underline w \in \R^{N_0 + N_L}$ of $w_h \in S_h^0(\Sigma).$ Here, 
for an easier implementation, we approximate the right-hand side 
$g \in H^{1/2}_{0,}(\Sigma)$ by $Q_h g$, where $Q_h$ is the $L^2$ projection on 
the space of piecewise linear, continuous functions fulfilling homogeneous 
initial conditions for $t=0.$ The assembling of the matrix 
$V_h \in \R^{(N_0 + N_L) \times (N_0 + N_L)}$ and the right-hand side 
$\underline g \in \R^{N_0 + N_L}$, i.e., the realization of $\mathcal H_T$, is done as proposed in
\cite[Subsection~2.2]{Zank2021Exact}. The integrals for computing the
projection $Q_h g$ are calculated by using high-order quadrature rules. 
The global linear system \eqref{Num:LGS} is solved by a direct solver.

In the numerical examples,
we consider the spatial domain $\Omega=(0,3)$, i.e., $L=3$, and
the time interval $(0,6)$, i.e., $T=6$. The lateral boundaries
$\{ 0 \} \times (0,T)$ and $ \{L \} \times (0,T)$ are discretized
uniformly into $N_0=N_L=2^{\ell+1}$ boundary elements each, $\ell=3,4,5,\ldots,12$.

In the first example, we consider the smooth solution
\[
  u_1(x,t) = \left \{ \begin{array}{lcl}
                        \frac{1}{2} \, (t-x-2)^3(x-t)^3
                        & & \mbox{for} \; x \leq t \leq 2 + x, \\[1mm]
                        0 & & \mbox{otherwise} .
                        \end{array} \right.
\]
Due to $w_1 = \partial_{n_x}u_1 \in H^1(\Sigma)$ and using the error estimate
\eqref{error L2}, we expect a linear order of convergence, as
confirmed by the numerical results given in Table
\ref{Table smooth}.

\begin{table}[h]

  \begin{center}
  \begin{tabular}{rrcc}
    \hline
    $\ell$ & $N_0+N_L$ & $\| w_1 - w_{1,h} \|_{L^2(\Sigma)}$ & eoc  \\
    \hline
 3 &     32 & 4.48 --1 &      \\ 
 4 &     64 & 2.11 --1 & 1.09  \\ 
 5 &    128 & 1.04 --1 & 1.02  \\ 
 6 &    256 & 5.18 --2 & 1.01  \\ 
 7 &    512 & 2.59 --2 & 1.00  \\ 
 8 &   1024 & 1.29 --2 & 1.00  \\ 
 9 &   2048 & 6.47 --3 & 1.00  \\ 
10 &   4096 & 3.23 --3 & 1.00  \\ 
11 &   8192 & 1.62 --3 & 1.00  \\ 
12 &  16384 & 8.09 --4 & 1.00  \\ 
    \hline
  \end{tabular}
\end{center}
  \caption{Numerical results for the boundary element method \eqref{Num:direktBEM} in the case
    $w_1 \in H^1(\Sigma)$.}
  \label{Table smooth}
\end{table}

\noindent
As a second example, we consider the singular solution
\[
  u_2(x,t) = \left \{ \begin{array}{ccl}
                        \displaystyle \frac{1}{2} \, | \sin (\pi(x-t))|
                        & & \mbox{for} \; x \leq t, \\[1mm]
                        0 & & \mbox{otherwise} ,
                        \end{array} \right .
\]
where we have $w_2 \in H^s(\Sigma)$ for $s < \frac{1}{2}$. Hence, 
using (\ref{error L2}), we expect the reduced order $\frac{1}{2}$ of
convergence when considering the error in $L^2(\Sigma)$. This is confirmed
by the numerical results as given in Table \ref{Table singular}.

\begin{table}[h]

  \begin{center}
  \begin{tabular}{rrcc}
    \hline
    $\ell$ & $N_0+N_L$ & $\| w_2 - w_{2,h} \|_{L^2(\Sigma)}$ & eoc  \\
    \hline
 3 &     32 & 2.59 +0 & 0.34  \\ 
 4 &     64 & 1.75 +0 & 0.56  \\ 
 5 &    128 & 1.21 +0 & 0.53  \\ 
 6 &    256 & 8.45 --1 & 0.52  \\ 
 7 &    512 & 5.93 --1 & 0.51  \\ 
 8 &   1024 & 4.18 --1 & 0.51  \\ 
 9 &   2048 & 2.95 --1 & 0.50  \\ 
10 &   4096 & 2.08 --1 & 0.50  \\ 
11 &   8192 & 1.47 --1 & 0.50  \\ 
12 &  16384 & 1.04 --1 & 0.50  \\     
    \hline
  \end{tabular}
\end{center}
  \caption{Numerical results for the boundary element method \eqref{Num:direktBEM} in the case
     $w_2 \in H^s(\Sigma)$, $s<\frac{1}{2}$.}
  \label{Table singular}
\end{table}

\section{Conclusions}
In this note, we have shown that the single layer boundary integral
operator of the wave equation in one space dimension is elliptic
in the energy space $[H^{1/2}_{,0}(\Sigma)]'$, when composed with
some modified Hilbert transformation. This result corresponds
to the well-known ellipticity results for boundary integral operators
related to second-order elliptic partial differential equations.
While this particular result is at this time restricted to the
spatially one-dimensional case, in the general case we were
already able to establish a related inf-sup stability condition
\cite{SteinbachUrzua:2021} instead. Although this is already
sufficient to do a numerical analysis of related boundary element methods,
it remains open whether we can prove ellipticity also in the
multi-dimensional case. It is obvious that we can extend this approach
also to the hypersingular boundary integral operator, and to the
double layer boundary integral operator. Ellipticity of boundary integral
operators is an important ingredient in the a priori and a posteriori
error analysis of boundary element methods, in the construction of
appropriate preconditioners, and in the coupling with finite element
methods. It goes without saying that this proposed new approach requires
more work in the numerical analysis, and in the implementation of
the proposed scheme, including the composition of the single layer
boundary integral operator and the modified Hilbert transformation, which are
both non-local. Nevertheless, this work may give some more insight into
the numerical analysis of existing boundary element methods for the
wave equation, and it presents an alternative approach for a reliable
and efficient numerical solution of the wave equation.

\section*{Acknowledgment}
The research of the second author was funded by the
John Fell Oxford University Press Research Fund.


\begin{thebibliography}{99}
\bibitem{Aimi:2009}
  A.~Aimi, M.~Diligenti, C.~Guardasoni, I.~Mazzieri, S.~Panizzi:
  An energy approach to space-time Galerkin BEM for wave propagation
  problems. Internat. J. Numer. Methods Engrg. 80 (2009) 1196--1240.
\bibitem{Bamberger:1986}
  A.~Bamberger, T.~Ha~Duong: Formulation variationnelle pour le calcul
  de la diffraction d’une onde acoustique par une surface rigide.
  Math. Meth. Appl. Sci. 8 (1986) 598--608.
\bibitem{ButzerTrebels:1968}
  P.~L.~Butzer, W.~Trebels: Hilberttransformation, gebrochene Integration
  und Differentiation, Springer Fachmedien Wiesbaden GmbH, 1968.
\bibitem{Costabel:2017}
  M.~Costabel, F.-J.~Sayas: Time-dependent problems with the boundary
  integral equation method. In: Encyclopedia of Computational Mechanics
  (E.~Stein, R.~Borst, T.~J.~R.~Hughes eds.), 2nd ed., Wiley, 2017.
\bibitem{Gimperlein:2017}
  H.~Gimperlein, Z.~Nezhi, E.~P.~Stephan: A priori error estimates for a
  time-dependent boundary element method for the acoustic wave equation
  in a half-space. Math. Methods Appl. Sci. 40 (2017) 448--462.
\bibitem{Hardy:1952}
  G.~Hardy, J.~E.~Littlewood, G.~P\'olya:
  Inequalities. Cambridge University Press, 1952.
\bibitem{Hassel:2017}
  M.~E.~Hassell, T.~Qiu, T.~S\'anchez-Vizuet, F.-J. Sayas: A new and
  improved analysis of the time domain boundary integral operators for
  the acoustic wave equation. J. Integral Equations Appl. 29 (2017) 107--136.
\bibitem{Joly:2017}
  P.~Joly, J.~Rodr\'iguez: Mathematical aspects of variational boundary
  integral equations for time dependent wave propagation.
  J. Integral Equations Appl. 29 (2017) 137--187.
\bibitem{Poelz:2021}
  D.~P\"olz, M.~Schanz: On the space-time discretization of variational
  retarded potential boundary integral equations. arXiv, 2103.16841v1, 2021.
\bibitem{Sayas:2013}
  F.-J.~Sayas: Energy estimates for Galerkin semidiscretizations of time
  domain boundary integral equations. Numer. Math. 124 (2013) 121--149.
\bibitem{Sayas:2016}
  F.-J.~Sayas: Retarded potentials and time domain boundary integral
  equations. A road map, volume 50 of Springer Series in Computational
  Mathematics. Springer, Cham, 2016.  
\bibitem{Steinbach:2008}
  O.~Steinbach: Numerical approximation methods for elliptic boundary
  value problems. Finite and boundary elements. Springer, New York, 2008. 
\bibitem{SteinbachUrzua:2021}
  O.~Steinbach, C.~Urz\'ua--Torres:
  A new approach to space-time boundary integral equations
  for the wave equation. arXiv, 2105.06800, 2021.
\bibitem{SteinbachZank:2016}
  O.~Steinbach, M.~Zank: Adaptive space-time boundary element methods
  for the wave equation. Proc. Appl. Math. Mech. 16 (2016) 777--778. 
\bibitem{SteinbachZank:2018}
  O.~Steinbach, M.~Zank: Coercive space--time finite element methods
  for initial boundary value problems.
  Electron. Trans. Numer. Anal. 52 (2020) 154--194.
\bibitem{SteinbachZank:2021}
  O.~Steinbach, M.~Zank: A note on the efficient evaluation of a
  modified Hilbert transformation. J. Numer. Math. 29 (2021) 47--61.
\bibitem{SteinbachZank:2021a}
  O.~Steinbach, M.~Zank: A generalized inf-sup stable variational
  formulation for the wave equation. arXiv, 2101.06293v1, 2021.
\bibitem{Urzua:2021}
  C.~Urz\'ua--Torres (joint work with O.~Steinbach):
  A new approach to time-domain boundary integral equations for the
  wave equation. Oberwolfach Reports 17 (2021) 371--373.
\bibitem{Zank2021Exact}
  M.~Zank: An exact realization of a modified {H}ilbert transformation for 
  space-time methods for parabolic evolution equations.
  Comput. Meth. Appl. Math. 21(2) (2021) 479--496.
\end{thebibliography}
\end{document}